\documentclass [twoside,reqno, 12pt] {amsart}

\usepackage{hyperref}

\usepackage{amsfonts}
\usepackage{amssymb}
\usepackage{a4}
\usepackage{color}
\usepackage{amscd,amssymb,amsfonts}

\usepackage{picture,xcolor}

\usepackage{etoolbox}

\usepackage{tikz}

\newtheorem{thm}{Theorem}[section]
\newtheorem{cor}[thm]{Corollary}

\newtheorem{rem}[thm]{Remark}

\theoremstyle{definition}

\numberwithin{equation}{section}

\newcommand{\R}{\mathbb{R}}

\def\tilde{\widetilde}
\def \bfo {\begin {eqnarray*} }
\def \efo {\end {eqnarray*} }
\def \ba {\begin {eqnarray*} }
\def \ea {\end {eqnarray*} }
\def \beq {\begin {eqnarray}}
\def \eeq {\end {eqnarray}}

\def \p {\partial}

\def\tilde{\widetilde}
\def \bfo {\begin {eqnarray*} }
\def \efo {\end {eqnarray*} }
\def \ba {\begin {eqnarray*} }
\def \ea {\end {eqnarray*} }
\def \beq {\begin {eqnarray}}
\def \eeq {\end {eqnarray}}

\def \p {\partial}

\begin{document}

\title[Fractional anisotropic Calder\'on problem on closed manifolds]{Fractional anisotropic Calder\'on problem on closed Riemannian manifolds}

\author[Feizmohammadi]{Ali Feizmohammadi}

\address{A. Feizmohammadi, The Fields Institute for Research in Mathematical Sciences\\
Toronto, Ontario M5T 3J1\\
Canada}

\email{afeizmoh@fields.utoronto.ca}

\author[Ghosh]{Tuhin Ghosh}
\address
        {T. Ghosh, Faculty of Mathematics\\
Bielefeld University, Postfach 100131\\
D-33501 Bielefeld, Germany}

\email{tghosh@math.uni-bielefeld.de}

\author[Krupchyk]{Katya Krupchyk}
\address
        {K. Krupchyk, Department of Mathematics\\
University of California, Irvine\\
CA 92697-3875, USA }

\email{katya.krupchyk@uci.edu}

\author[Uhlmann]{Gunther Uhlmann}

\address
       {G. Uhlmann, Department of Mathematics\\
       University of Washington\\
       Seattle, WA  98195-4350\\
       USA\\
        and Institute for Advanced Study of the Hong Kong University of Science and Technology}
\email{gunther@math.washington.edu}

\maketitle

\begin{abstract}
In this paper we solve the fractional anisotropic Calder\'on problem on closed Riemannian manifolds of dimensions two and higher. Specifically, we prove that the knowledge of the local source--to--solution map for the fractional Laplacian, given on an arbitrary small open nonempty a priori known subset of a smooth closed connected Riemannian manifold, determines the Riemannian manifold up to an isometry.  This can be viewed as a nonlocal analog of the anisotropic Calder\'on problem in the setting of closed Riemannian manifolds, which is wide open in dimensions three and higher.  
\end{abstract}

\section{Introduction and statement of results}

The anisotropic Calder\'on problem consists of recovering the electrical conductivity matrix of a medium, up to a change of variables, from the current and voltage measurements performed along the boundary of the medium, see \cite{Calderon_1980}, \cite{Uhlmann_2014}.  More generally, the problem can be stated on a compact Riemannian manifold with boundary as follows. Let $(M, g)$ be a smooth connected compact Riemannian manifold of dimension $n\ge 2$ with smooth boundary $\p M$, and let $-\Delta_g$ be the positive Laplace--Beltrami operator on $M$. Consider the Cauchy data on $\partial M$ of harmonic functions on $M$, 
\[
C_g:=\{(u|_{\p M}, \p_\nu u|_{\p M}): u\in C^\infty(M) \text{ such that } -\Delta_g u=0 \text{ in } M\}.
\]
If $\psi:M\to M$ is a diffeomorphism such that $\psi|_{\p M}=Id$ then $C_{\psi^*g}=C_g$, see \cite{Kohn_Vogelius_1983}, \cite{Lee_Uhlmann_1989}. The geometric version of the anisotropic Calder\'on problem asks to determine a smooth connected compact Riemannian manifold $(M,g)$ uniquely, up to such an isometry. While in dimension $n=2$, the geometric version of the anisotropic Calder\'on problem, with an additional obstruction arising from the conformal invariance of the Laplacian, is settled in  \cite{Lassas_Uhlmann_2001}, in dimensions $n\ge 3$ it is wide open and represents one of the major open problems in the field of inverse problems.  In dimensions $n\ge 3$, it is solved for real-analytic manifolds in \cite{Lassas_Taylor_Uhlmann_2003}, \cite{Lassas_Uhlmann_2001}, \cite{Lee_Uhlmann_1989},  see also  \cite{Guillarmou_Sa_Barreto_2009}.  In the $C^\infty$ case, the uniqueness for the anisotropic Calder\'on problem has only been established in \cite{DKSaloU_2009}, \cite{DKurylevLS_2016} for metrics in a fixed conformal class in the case of a conformally transversally anisotropic manifold and under the assumption that the geodesic ray transform on the transversal manifold is injective.

The purpose of this paper is to provide a solution of the fractional anisotropic Calder\'on problem on a closed smooth connected Riemannian manifold. Here in contrast to the anisotropic Calder\'on problem, we are able to recover uniquely a smooth closed connected Riemannian manifold up to a natural obstruction in any dimension $\ge 2$. By a closed manifold we mean a compact manifold without boundary. To state the problem in precise terms, let $(M,g)$ be a smooth closed connected Riemannian manifold of dimension $n\ge 2$. Let $-\Delta_g$ be the positive Laplace--Beltrami operator on $M$. It is a self-adjoint operator on $L^2(M)$ with the domain $\mathcal{D}(-\Delta_g)=H^2(M)$, the standard Sobolev space on $M$. We denote by $0=\lambda_0<\lambda_1<\lambda_2<\dots$ the distinct eigenvalues of $-\Delta_g$ and  by $d_k$ the multiplicity of $\lambda_k$, $k=0,1,2,\dots$.  Let $\varphi_{k,1},\dots, \varphi_{k,d_k}$ be an $L^2(M)$--orthonormal basis for the eigenspace $\text{Ker}(-\Delta_g-\lambda_k)$ corresponding to $\lambda_k$, $k=0,1,2,\dots$, and let $\pi_k:L^2(M)\to \text{Ker}(-\Delta_g-\lambda_k)$ be the orthogonal projection onto the eigenspace of $\lambda_k$ given by 
\[
\pi_k f=\sum_{j=1}^{d_k}(f,\varphi_{k,j})_{L^2(M)}\varphi_{k,j}, \quad f\in L^2(M),
\]
$k=0,1,2,\dots$. Here $(\cdot, \cdot)_{L^2(M)}$ is the $L^2$ inner product on $M$. 

Let $\alpha\in (0,1)$. By the spectral theorem we define the fractional Laplacian $(-\Delta_g)^\alpha$ of order $\alpha$ as an unbounded self-adjoint operator on $L^2(M)$ given by
\[
(-\Delta_g)^\alpha u=\sum_{k=0}^\infty \lambda_k^\alpha \pi_k u,
\]
equipped with the domain $\mathcal{D}((-\Delta_g)^\alpha):=\{u\in L^2(M):\sum_{k=0}^\infty \lambda_k^{2\alpha}\|\pi_k u\|_{L^2(M)}^2<\infty\}=H^{2\alpha}(M)$.

Let $\mathcal{O}\subset M$ be an open nonempty set and let $f\in C^\infty_0(\mathcal{O})$ be such that 
\[
(f,1)_{L^2(M)}=0.
\]
 Then the equation 
\begin{equation}
\label{eq_int_2}
(-\Delta_g)^\alpha u=f\quad \text{in}\quad M
\end{equation}
has a unique solution $u=u^f\in C^\infty(M)$ with the property that $(u^f, 1)_{L^2(M)}=0$, given by 
\[
u^f=(-\Delta_g)^{-\alpha}f=\sum_{k=1}^\infty \lambda_k^{-\alpha} \pi_k f.
\]
Associated to \eqref{eq_int_2} we define the local source--to--solution map $L_{M,g,\mathcal{O}}$ by 
\begin{equation}
\label{eq_int_3}
L_{M,g,\mathcal{O}}(f):=u^f|_{\mathcal{O}}=((-\Delta_g)^{-\alpha}f)|_{\mathcal{O}}, 
\end{equation}
where $f\in C^\infty_0(\mathcal{O})$ with $(f,1)_{L^2(M)}=0$ and $u^f\in C^\infty(M)$ is the unique solution to \eqref{eq_int_2}  with the property that $(u^f, 1)_{L^2(M)}=0$. 

The fractional anisotropic Calder\'on inverse problem that we are interested in is as follows: does the knowledge of the local source--to--solution map $L_{M,g,\mathcal{O}}$, the observation set $\mathcal{O}$, and the metric $g|_{\mathcal{O}}$  determine the manifold $(M,g)$ globally? It is known that there is an obstruction to uniqueness that comes from isometries on $M$ that fix the set $\mathcal{O}$, see \cite[Theorem 4.2]{Ghosh_Uhlmann_2021}.  Thus, the best that one could expect is to determine the manifold up to this obstruction. Recently, it was established in \cite{Feizmohammadi_2021} that one can indeed recover the isometric class of the manifold $(M,g)$ assuming that the observation set $(\mathcal{O}, g|_{\mathcal{O}})$ locally admits a metric in Gevrey class with some Gevrey index belonging to $[1,2)$.  

In this paper we provide a complete solution to the fractional anisotropic Calder\'on problem on a smooth closed connected Riemannian manifold. Our main result is as follows.
\begin{thm}
\label{thm_main}
Let $\alpha\in (0,1)$. Let $(M_1, g_1)$ and $(M_2,g_2)$  be smooth closed connected Riemannian manifolds of dimension $n\ge 2$, and let  $\mathcal{O}_j\subset M_j$, $j=1,2$, be open nonempty sets. Assume that 
\begin{equation}
\label{eq_int_4}
(\mathcal{O}_1, g_1|_{\mathcal{O}_1})=(\mathcal{O}_2, g_2|_{\mathcal{O}_2}):=(\mathcal{O}, g).
\end{equation}
Assume furthermore that 
\begin{equation}
\label{eq_int_5}
L_{M_2,g_2,\mathcal{O}_2}(f)= L_{M_1,g_1,\mathcal{O}_1}( f),
\end{equation}
for all $f\in C^\infty_0(\mathcal{O})$ with $(f,1)_{L^2(\mathcal{O})}=0$. Then there exists a diffeomorphism $\Phi:M_1\to M_2$  such that $\Phi^*g_2=g_1$ on $M_1$. 
\end{thm}

Theorem \ref{thm_main}  is of main interest when the set $\mathcal{O}$ is sufficiently small and in this case, we have the following consequence of Theorem \ref{thm_main}. 
\begin{cor}
\label{cor_main}
Under the assumptions of Theorem \ref{thm_main} and assuming that  $\mathcal{O}$ is sufficiently small, there exists a diffeomorphism $\Phi:M_1\to M_2$  such that $\Phi^*g_2=g_1$ on $M_1$ and $\Phi(x)=x$ for $x\in \mathcal{O}$. 
\end{cor}

\begin{rem}
In two--dimensional case, while there is an additional obstruction in the geometric version of the anisotropic Calder\'on problem, coming from the conformal invariance of the Laplacian, see \cite{Lassas_Uhlmann_2001}, this obstruction is not present in the fractional anisotropic Calderon problem, see Theorem \ref{thm_main} and Corollary \ref{cor_main}.
\end{rem}

The study of the fractional Calder\'on problem was initiated in \cite{Ghosh_Salo_Uhlmann_2020} where the unknown potential in the fractional Schr\"odinger equation on a bounded domain in the Euclidean space was determined from exterior measurements. Following this work, inverse problems of recovering lower order terms for fractional elliptic equations have been studied extensively, see for example \cite{Ghosh_Ruland_Salo_Uhlmann_2020}, \cite{Ghosh_Lin_Xiao_2017}, \cite{Ruland_Salo_2020}, \cite{Ruland_Salo_2018}, \cite{Bhattacharyya_Ghosh_Uhlmann_2021}, \cite{Cekic_Lin_Ruland_2020}, \cite{Covi_2020}, \cite{Covi_2020_2}, \cite{Covi_preprint}, \cite{Covi_Garcia-Ferrero_Ruland_preprint},  \cite{Covi_M_R_Uhlmann_preprint}, \cite{Li_Li_2020}, \cite{Li_Li_2021},  \cite{Ruland_2021} for some of the important contributions. In all of those works, it is assumed that the leading order coefficients are known.

\begin{rem}
To the best of our knowledge, the inverse problem of recovering the anisotropic leading order coefficients appearing in the fractional powers of self-adjoint positive operators has only been studied in the recent works  \cite{Ghosh_Uhlmann_2021} and \cite{Feizmohammadi_2021}.  The work \cite{Ghosh_Uhlmann_2021} deals with the inverse problem of determining  anisotropic leading order coefficients arising in the fractional powers of general self-adjoint positive second order operators on a bounded domain in the Euclidean space from exterior measurements.  Specifically, it was proved in  \cite{Ghosh_Uhlmann_2021} that the exterior data for the nonlocal inverse problem determines the boundary Cauchy data of solutions to the corresponding local equation, and using this determination, the results available for the local anisotropic inverse problems were extended to the case of nonlocal anisotropic inverse problems.  The recent work \cite{Feizmohammadi_2021} solves the fractional anisotropic Calder\'on problem on a smooth closed connected Riemannian manifold under the assumption that the observation set $(\mathcal{O}, g|_{\mathcal{O}})$ locally admits a metric in Gevrey class with some Gevrey index belonging to $[1,2)$.  Our Theorem \ref{thm_main} and Corollary \ref{cor_main} can be viewed as a significant extension of the latter result. 
\end{rem}

Let us also remark that geometric inverse problems on closed Riemannian manifolds with measurements performed on an open subset of the manifold were studied in \cite{Krupchyk_Kurylev_Lassas_2008},  \cite{HLOS_2018}, \cite{HLYZ_2020}, \cite{Bosi_Kurylev_Lassas}, \cite{Kurylev_Lassas_Yamaguchi}, \cite{FILLN} in the context of inverse interior spectral problems, as well as inverse problems for the wave and heat equations. 

Let us now proceed to describe the main ideas in the proof of Theorem \ref{thm_main}. To that end, let us write $P_{g_j}:=-\Delta_{g_j}$, $j=1,2$, and write $e^{-tP_g}$, $t\ge 0$, for the heat semigroup, and $e^{-tP_g}(\cdot,\cdot)$ for its Schwartz kernel. We have $e^{-tP_g}(x,y)\in C^\infty((0,\infty)\times M\times M)$, see  \cite[Theorem 3.5]{Strichartz_1983}. 
First we pass from the equality of the local source--to-solutions maps \eqref{eq_int_5} on the open set $\mathcal{O}$ to the equality for the heat kernels of the operators $P_{g_j}$ on $\mathcal{O}$, 
\begin{equation}
\label{eq_int_6}
e^{-tP_{g_1}}(x,y)=e^{-tP_{g_2}}(x,y), \quad x,y\in \mathcal{O}, \quad t>0.
\end{equation}
In doing so, we use a representation of the operator $P_{g_j}^{-\alpha}$, defined on $L^2(M_j)$--functions which are orthogonal to the constant function $1$,  in terms of the heat semigroup $e^{-t P_{g_j}}$. Theorem \ref{thm_main} is therefore a consequence of the following result, which we state in the more general setting of complete Riemannian manifolds.  To that end, let $(N,g)$ be a smooth complete Riemannian manifold of dimension $n\ge 2$ without boundary. Since $(N,g)$ is complete, the Laplace--Beltrami operator $P_g:=-\Delta_g$ associated with the metric $g$, is a nonnegative self-adjoint operator on $L^2(N)$, when equipped with the maximal domain, see  \cite[Theorem 2.4]{Strichartz_1983}. We have $e^{-tP_g}(x,y)\in C^\infty((0,\infty)\times N\times N)$, see  \cite[Theorem 3.5]{Strichartz_1983}. We have the following result. 

\begin{thm}
\label{thm_from_heat_semigroup_to_wave}
Let $(N_1,g_1)$ and $(N_2,g_2)$ be  smooth connected complete Riemannian manifolds of dimension $n\ge 2$ without boundary. Let $\mathcal{O}_j\subset N_j$, $j=1,2$, be open nonempty sets. Assume that $\mathcal{O}_1=\mathcal{O}_2:=\mathcal{O}$. Assume furthermore that 
\begin{equation}
\label{eq_5_2}
e^{-tP_{g_1}}(x,y)=e^{-tP_{g_2}}(x,y), \quad x,y\in \mathcal{O}, \quad t>0.
\end{equation}
Then there exists a diffeomorphism $\Phi:N_1\to N_2$  such that $\Phi^*g_2=g_1$ on $N_1$. 
\end{thm}

Theorem \ref{thm_from_heat_semigroup_to_wave} is known in the case of closed Riemannian manifolds and it follows from \cite{Krupchyk_Kurylev_Lassas_2008},  \cite{HLOS_2018}, \cite{HLYZ_2020}.  To the best of our knowledge, Theorem 
\ref{thm_from_heat_semigroup_to_wave} is new in the complete case.  We present a proof of Theorem \ref{thm_from_heat_semigroup_to_wave}, which is different from those available in the compact case,  for completeness and convenience of the reader. Furthermore, we hope that Theorem \ref{thm_from_heat_semigroup_to_wave} could be useful for some other geometric inverse problems in the more general setting of complete Riemannian manifolds.  To prove Theorem 
\ref{thm_from_heat_semigroup_to_wave}, we shall reduce the problem to an inverse problem for the wave equation with interior measurements on $\mathcal{O}$, and apply \cite[Theorem 2]{HLOS_2018} to conclude that the manifolds $(N_1,g_1)$ and $(N_2,g_2)$ are Riemannian isometric. When performing the reduction working on a complete manifold, we rely on the functional calculus of self-adjoint operators, making use, in particular, of a variant of the transmutation formula of Kannai, see \cite{Kannai_1977}.  We refer to \cite[Sections 9.2, 9.5]{Isakov_book} for the use of the transmutation formula of \cite{Kannai_1977} in inverse problems.

The plan of this paper is as follows. Section \ref{sec_proof_main_result} is devoted to the proof of Theorem \ref{thm_main} and Corollary \ref{cor_main}. Section \ref{sec_proof_from_heat_to_wave} contains the proof of Theorem \ref{thm_from_heat_semigroup_to_wave}.

\section{Proof of Theorem \ref{thm_main} and Corollary \ref{cor_main}}

\label{sec_proof_main_result}

\subsection{Proof of Theorem \ref{thm_main} }
Let $P_{g_j}=-\Delta_{g_j}$, $j=1,2$,  and let $\alpha\in (0,1)$. Let us write 
\[
P_{g_j}^{-\alpha} v_=\sum_{k=1}^\infty (\lambda_k^{(j)})^{-\alpha}\pi_k^{(j)} v_j,
\]
where $v_j\in L^2(M_j)$ is such that $(v_j,1)_{L^2(M_j)}=0$, $j=1,2$. It will be convenient for us to express the self-adjoint operator $P_{g_j}^{-\alpha}$ in terms of the heat semigroup $e^{-tP_{g_j}}$, $t\ge 0$, when acting on the orthogonal complement of the one-dimensional subspace of $L^2(M_j)$, spanned by the constant function $1$.  To that end, using the definition of the Gamma function 
\[
\Gamma(\alpha)=\int_0^\infty e^{-t}t^{\alpha-1}dt, 
\]
we get 
\begin{equation}
\label{eq_2_1}
a^{-\alpha}=\frac{1}{\Gamma(\alpha)}\int_0^\infty e^{-at}\frac{1}{t^{1-\alpha}}dt, \quad a>0.
\end{equation}
It follows from \eqref{eq_2_1} that 
\begin{equation}
\label{eq_2_2}
P_{g_j}^{-\alpha}v_j=\frac{1}{\Gamma(\alpha)}\int_0^\infty e^{-tP_{g_j}}v_j\frac{1}{t^{1-\alpha}}dt,
\end{equation}
where $v_j\in L^2(M_j)$ is such that $(v_j,1)_{L^2(M_j)}=0$. The integral in \eqref{eq_2_2} converges in $L^2(M_j)$. 

Let $\omega_1\subset\subset \mathcal{O}$ be an open nonempty set. Then there exists an open nonempty set $\omega_2\subset \mathcal{O}$ such that $\overline{\omega_1}\cap \overline{\omega_2}=\emptyset$. Let  $f\in C^\infty_0(\omega_1)$. By \eqref{eq_int_4}, for $m=1,2,\dots$,  we have
\[
\Delta_{g_1}^mf=\Delta_{g_2}^mf=\Delta_g^m f \quad \text{on}\quad \omega_1. 
\] 
Furthermore, for $m=1,2,\dots$, we have that the functions $\Delta_g^m f \in C^\infty_0(\omega_1)$ are such that  $(\Delta_g^m f, 1)_{L^2(M_j)}=0$, $j=1,2$. Thus, it follows from 
\eqref{eq_int_5} and \eqref{eq_int_3} that for $m=1,2,\dots$, 
\begin{equation}
\label{eq_2_3}
\big(P_{g_1}^{-\alpha}\Delta_g^m f \big)|_{\mathcal{O}}=\big(P_{g_2}^{-\alpha}\Delta_g^m f \big)|_{\mathcal{O}}.
\end{equation}
Using \eqref{eq_2_2}, we get from \eqref {eq_2_3} that 
\begin{equation}
\label{eq_2_6}
\int_0^\infty \big((e^{-tP_{g_1}}-e^{-tP_{g_2}})\Delta_g^m f\big)(x)\frac{dt}{t^{1-\alpha}}=0,
\end{equation}
for $x\in \mathcal{O}$, and $m=1,2,\dots$.  

Next we shall argue as in \cite[Proof of Proposition 3.1]{Ghosh_Uhlmann_2021}. Using that the function $t\mapsto e^{-tP_{g_j}} f\in C^\infty([0,\infty), L^2(M_j))$, and that $e^{-tP_{g_j}}\Delta_{g_j}^m= \Delta_{g_j}^m e^{-tP_{g_j}}$ for all $t\ge 0$ on $\mathcal{D}(\Delta_{g_j}^m)$,  we get for any $m=1,2,\dots$, 
\begin{equation}
\label{eq_2_6_1}
\big(e^{-tP_{g_j}}\Delta_g^m f\big)(x)=\p_t^m \big(e^{-tP_{g_j}} f\big)(x), 
\end{equation}
for  $x\in \mathcal{O}$. Combining \eqref{eq_2_6_1} and  \eqref{eq_2_6}, we obtain  that 
\begin{equation}
\label{eq_2_7}
\int_0^\infty \p_t^m \big((e^{-tP_{g_1}}-e^{-tP_{g_2}})f\big)(x)\frac{dt}{t^{1-\alpha}}=0,
\end{equation}
for  $x\in \mathcal{O}$ and $m=1,2,\dots$.  Assuming that $x\in \omega_2$, we shall integrate by parts in \eqref{eq_2_7} $m$ times. We claim that there will be no contributions from the end points. Indeed, for $l=0,\dots, m-1$, using \eqref{eq_2_6_1}, we have for $t>0$ and $x\in \omega_2$,
\begin{equation}
\label{eq_2_8}
\begin{aligned}
\p_t^{l} \big((e^{-tP_{g_1}}-e^{-tP_{g_2}})f\big)(x)=\big((e^{-tP_{g_1}}-e^{-tP_{g_2}})\Delta_g^{l} f\big)(x)\\
=\int_{\omega_1} (e^{-tP_{g_1}}(x,y)-e^{-tP_{g_2}}(x,y))(\Delta_g^{l}f )(y)dV_g(y),
\end{aligned}
\end{equation}
where $e^{-tP_{g_j}}(x,y)$ is the kernel of $e^{-tP_{g_j}}$, $j=1,2$, and $dV_g$ is the Riemannian volume element. It follows from \eqref{eq_2_8} that for $t>0$ and $x\in \omega_2$, 
\begin{equation}
\label{eq_2_9}
\big|\p_t^{l} \big((e^{-tP_{g_1}}-e^{-tP_{g_2}})f\big)(x) \big|\le \|e^{-tP_{g_1}}(\cdot,\cdot)-e^{-tP_{g_2}}(\cdot,\cdot)\|_{L^\infty(\omega_2\times \omega_1)}\|\Delta_g^{l}f \|_{L^1(\omega_1)},
\end{equation}
$l=0,\dots, m-1$. 
To proceed we need the following  pointwise upper Gaussian estimate on the heat kernel
\begin{equation}
\label{eq_2_9_1}
|e^{-tP_{g_j}}(x,y)|\le Ct^{-\frac{n}{2}}e^{-\frac{cd_{g_j}(x,y)^2}{t}}, \quad t>0, \quad x,y\in M_j,
\end{equation}
which follows from 
\cite{Grigoryan_1997}, combining with the estimate of \cite{Varopoulos_1985},
\[
\|e^{-tP_{g_j}}\|_{L^1(M_j)\to L^\infty(M_j)}\le Ct^{-\frac{n}{2}}, \quad t>0.
\]
Here  $c>0$ and $d_{g_j}(\cdot, \cdot)$ denotes the Riemannian distance on $(M_j,g_j)$, $j=1,2$.  Using \eqref{eq_2_9_1}, we get from \eqref{eq_2_9} that for $0<t<1$ and $x\in \omega_2$, 
\begin{equation}
\label{eq_2_10}
\big|\p_t^{l} \big((e^{-tP_{g_1}}-e^{-tP_{g_2}})f\big)(x) \big| \le Ce^{-\frac{\tilde c}{t}}\|\Delta_g^{l}f \|_{L^1(\omega_1)},
\end{equation}
where $\tilde c>0$ depends on $d_g(\overline{\omega_1}, \overline{\omega_2})>0$, $l=0,\dots, m-1$.  For $t>1$ and $x\in \omega_2$, 
we get 
\begin{equation}
\label{eq_2_10_1}
\big|\p_t^{l} \big((e^{-tP_{g_1}}-e^{-tP_{g_2}})f\big)(x) \big| \le Ct^{-\frac{n}{2}}\|\Delta_g^{l}f \|_{L^1(\omega_1)},
\end{equation}
$l=0,\dots, m-1$. 
The bounds \eqref{eq_2_10} and \eqref{eq_2_10_1} show the claim. Thus, integrating by parts $m$ times in \eqref{eq_2_7}, we obtain 
\begin{equation}
\label{eq_2_5}
\int_0^\infty \big((e^{-tP_{g_1}}-e^{-tP_{g_2}})f\big)(x)\frac{dt}{t^{1+m-\alpha}}=0,
\end{equation}
for  $x\in \omega_2$ and $m=1,2,\dots$. Rewriting \eqref{eq_2_5} as follows
\begin{equation}
\label{eq_2_5_1}
\int_0^\infty \big((e^{-tP_{g_1}}-e^{-tP_{g_2}})f\big)(x)\frac{dt}{t^{2+m-\alpha}}=0,
\end{equation}
for $x\in \omega_2$ and $m=0,1,2,\dots$, and making the change of variables $s=1/t$, we get 
\begin{equation}
\label{eq_2_5_2}
\int_0^\infty \varphi(s)s^{m}ds=0,
\end{equation}
where 
\[
\varphi(s)=\frac{\big((e^{-\frac{1}{s}P_{g_1}}-e^{-\frac{1}{s}P_{g_2}})f\big)(x)}{s^\alpha}, \quad x\in \omega_2.
\]
In view of \eqref{eq_2_10} and \eqref{eq_2_10_1}, we have for $s\ge 0$, 
\[
|\varphi(s)|\le \mathcal{O}(1)\frac{e^{-cs}}{s^\alpha}, 
\]
where $c>0$. It follows that the Fourier transform of $1_{[0,\infty)}\varphi$, 
\[
\mathcal{F}(1_{[0,\infty)}\varphi)(\xi) = \int_0^\infty \varphi(s)e^{-i\xi s}ds
\]
is holomorphic for $\text{Im}\ \xi<c$. In view of \eqref{eq_2_5_2}, $\mathcal{F}(1_{[0,\infty)}\varphi)$ vanishes at $0$ with all derivatives, and therefore, $\varphi(s)=0$ for $s>0$. Hence, we get 
\begin{equation}
\label{eq_2_12}
 \big((e^{-tP_{g_1}}-e^{-tP_{g_2}})f\big)(x)=0,
\end{equation}
for  $t>0$ and $x\in \omega_2$.  Furthermore, the function $\big((e^{-tP_{g_1}}-e^{-tP_{g_2}})f\big)\big|_{(0,\infty)\times \overline{\mathcal{O}}}\in C^\infty((0,\infty)\times \overline{\mathcal{O}})$ satisfies the heat equation 
\begin{equation}
\label{eq_2_13}
(\p_t-\Delta_g)((e^{-tP_{g_1}}-e^{-tP_{g_2}})f)=0\quad \text{in}\quad \mathcal{O}. 
\end{equation}
Assuming without loss of generality that $\mathcal{O}$ is connected and contained in a single coordinate patch for both manifolds $M_1$ and $M_2$,  in view of \eqref{eq_2_12} and \eqref{eq_2_13}, by the unique continuation for the heat equation, see \cite[Sections 1 and 4]{Lin_1990}, we obtain that 
\begin{equation}
\label{eq_2_14}
 \big((e^{-tP_{g_1}}-e^{-tP_{g_2}})f\big)(x)=0,
\end{equation}
for $t>0$ and  $x\in \mathcal{O}$.  
Recalling that here $f\in C^\infty_0(\omega_1)$ is arbitrary and $\omega_1\subset\subset \mathcal{O}$  is also arbitrary, we conclude from \eqref{eq_2_14} that
\begin{equation}
\label{eq_2_15}
 e^{-tP_{g_1}}f|_{\mathcal{O}}=e^{-tP_{g_2}}f|_{\mathcal{O}},
\end{equation}
for $t>0$ and $f\in C^\infty_0(\mathcal{O})$. Now \eqref{eq_2_15} implies that 
\[
e^{-tP_{g_1}}(x,y)=e^{-tP_{g_2}}(x,y), \quad x,y\in \mathcal{O}, \quad t>0,
\]
and therefore, an application of Theorem \ref{thm_from_heat_semigroup_to_wave} completes the proof of Theorem \ref{thm_main}. 

\subsection{Proof of Corollary \ref{cor_main}}
Corollary \ref{cor_main} follows from Theorem \ref{thm_main} as well as the observation that the proof of \cite[Theorem 2]{HLOS_2018} used in our proof of Theorem \ref{thm_from_heat_semigroup_to_wave} gives the existence of a diffeomorphism $\Phi:M_1\to M_2$  such that $\Phi^*g_2=g_1$ on $M_1$ and $\Phi(x)=x$ for $x\in \mathcal{O}$, provided that $\mathcal{O}$ is sufficiently small.

\section{Proof of Theorem \ref{thm_from_heat_semigroup_to_wave}}
\label{sec_proof_from_heat_to_wave}

We shall reduce the problem to an inverse problem for the wave equation with interior measurements on $\mathcal{O}$. To that end, let $F\in C^\infty_0((0,\infty)\times \mathcal{O})$ and consider the following inhomogeneous initial value problem for the wave equation,
\begin{equation}
\label{eq_5_3}
\begin{cases}
(\p_t^2-\Delta_{g_j})u_j(t,x)=F(t,x), & (t,x)\in (0,\infty)\times N_j,\\
u_j(0,x)=0, & x\in N_j,\\
\p_t u_j(0,x)=0, & x\in N_j,
\end{cases}
\end{equation}
$j=1,2$. The problem \eqref{eq_5_3} has a unique solution $u_j=u_j^F\in C^\infty([0,\infty)\times N_j)$, see \cite{HLOS_2018}. Associated to \eqref{eq_5_3}, we define the local source--to--solution map on $\mathcal{O}$ by
\begin{equation}
\label{eq_5_4}
L^{\text{wave}}_{N_j,g_j, \mathcal{O}}:C^\infty_0((0,\infty)\times \mathcal{O})\to C^\infty ([0,\infty)\times \mathcal{O}), \quad L^{\text{wave}}_{N_j,g_j, \mathcal{O}}(F)=u_j^F|_{\mathcal{O}},
\end{equation}
where  $u_j^F$ is the unique solution to \eqref{eq_5_3}.  Our goal is to show that the equality \eqref{eq_5_2} for the heat kernels implies the equality of the local source--to--solution maps, 
\begin{equation}
\label{eq_5_5}
L^{\text{wave}}_{N_1,g_1, \mathcal{O}}(F)= L^{\text{wave}}_{N_2,g_2, \mathcal{O}}(F),
\end{equation}
for $F\in C^\infty_0((0,\infty)\times \mathcal{O})$, and thus, the conclusion of Theorem \ref{thm_from_heat_semigroup_to_wave} will follow from \cite[Theorem 2]{HLOS_2018}. 

To prove \eqref{eq_5_5}, we shall use functional calculus of self-adjoint operators. Specifically, the solution to \eqref{eq_5_3} is given by 
\begin{equation}
\label{eq_5_6}
u_j^F(t,x)=\int_0^t \frac{\sin((t-s)\sqrt{P_{g_j}})}{\sqrt{P_{g_j}}}F(s,x)ds, \quad (t,x)\in [0,\infty)\times N_j,
\end{equation}
$j=1,2$.  To proceed we shall use a variant of the transmutation formula of Kannai, which transforms the solution to the wave equation to the solution of the heat equation, see \cite{Kannai_1977}. To obtain the needed transmutation formula, we first write the Fourier inversion formula for the Gaussian $\R\ni \lambda\mapsto e^{-t\lambda^2}$, $t>0$, 
\begin{equation}
\label{eq_5_7}
e^{-t\lambda^2}=\frac{1}{(4\pi t)^{1/2}}\int_{-\infty}^{+\infty} e^{-\frac{s^2}{4t}}e^{is\lambda}ds=
\frac{1}{(4\pi t)^{1/2}}\int_{-\infty}^{+\infty} e^{-\frac{s^2}{4t}}\cos(s\lambda)ds, \quad t>0.
\end{equation}
Integrating by parts in \eqref{eq_5_7} for $\lambda\ne 0$, we get 
\begin{equation}
\label{eq_5_8}
\begin{aligned}
e^{-t\lambda^2}=
\frac{1}{4\pi^{1/2}t^{3/2}}\int_{-\infty}^{+\infty} s e^{-\frac{s^2}{4t}}\frac{\sin(s\lambda)}{\lambda}ds=\frac{2}{4\pi^{1/2}t^{3/2}}\int_{0}^{+\infty} s e^{-\frac{s^2}{4t}}\frac{\sin(s\lambda)}{\lambda}ds\\
=\frac{1}{4\pi^{1/2}t^{3/2}}\int_{0}^{+\infty} e^{-\frac{\tau}{4t}}\frac{\sin(\sqrt{\tau}\lambda)}{\lambda}d\tau.
\end{aligned}
\end{equation}
Note that in the last equality in \eqref{eq_5_8} we have used the change of variables $\tau=s^2$. The formula \eqref{eq_5_8} holds also for $\lambda=0$ in view of the dominated convergence theorem.  Using \eqref{eq_5_8}, we have the following variant of the transmutation formula of Kannai, 
\begin{equation}
\label{eq_5_9}
e^{-tP_{g_j}}v_j=\frac{1}{4\pi^{1/2}t^{3/2}}\int_{0}^{+\infty} e^{-\frac{\tau}{4t}}\frac{\sin(\sqrt{\tau}\sqrt{P_{g_j}})}{\sqrt{P_{g_j}}}v_jd\tau,  \quad t>0,
\end{equation}
where $v_j\in L^2(N_j)$, $j=1,2$.

Let $f\in C^\infty_0(\mathcal{O})$. Then \eqref{eq_5_2} and \eqref{eq_5_9} imply that 
\begin{equation}
\label{eq_5_10}
\int_{0}^{+\infty} e^{-\tau t} \bigg(\frac{\sin(\sqrt{\tau}\sqrt{P_{g_1}})}{\sqrt{P_{g_1}}}f\bigg)(x)d\tau=\int_{0}^{+\infty} e^{-\tau t} \bigg(\frac{\sin(\sqrt{\tau}\sqrt{P_{g_2}})}{\sqrt{P_{g_2}}}f\bigg)(x)d\tau,  
\end{equation}
for $t>0$ and $x\in \mathcal{O}$. Inverting the Laplace transform in \eqref{eq_5_10}, we get 
\begin{equation}
\label{eq_5_11}
\bigg(\frac{\sin(t\sqrt{P_{g_1}})}{\sqrt{P_{g_1}}}f\bigg)(x)=\bigg(\frac{\sin(t\sqrt{P_{g_2}})}{\sqrt{P_{g_2}}}f\bigg)(x),
\end{equation}
for $t>0$ and $x\in \mathcal{O}$. 
Recalling that $F\in C^\infty_0((0,\infty)\times \mathcal{O})$, we conclude from \eqref{eq_5_11} and \eqref{eq_5_6} that 
\[
u_1^F(t,x)=u_2^F(t,x), 
\]
for $t>0$ and $x\in \mathcal{O}$, and thus, in view of \eqref{eq_5_4}, we get \eqref{eq_5_5}. This completes the proof of Theorem \ref{thm_from_heat_semigroup_to_wave}.

\section*{Acknowledgements}

A.F. gratefully acknowledges support from the Fields institute for research in mathematical sciences.  The research of T.G. is supported by the Collaborative Research Center, membership no. 1283, Universit\"at Bielefeld.  The research of K.K. is partially supported by the National Science Foundation (DMS 2109199). The research of G.U. is partially supported by NSF, a Walker Professorship at UW, a Si-Yuan Professorship at IAS, HKUST, and Simons Fellowship. Part of this research  was performed while G.U. was visiting the Institute for Pure and Applied Mathematics (IPAM), which is supported by the National Science Foundation (Grant No. DMS-1925919).

\end{document}